\documentclass[11pt,a4paper]{article}

\usepackage{amsmath}
\usepackage{amsthm}
\usepackage{amsfonts}

\parskip\smallskipamount
\newtheorem{theorem}{Theorem}[section]
\newtheorem{lem}{Lemma}[section]

\newtheorem{prop}{Proposition}[section]
\theoremstyle{remark}

\newtheorem{rem}{Remark}[section]
\newtheorem{example}{Example}[section]

\renewcommand{\Pr}{\mathbf{P}}
\newcommand{\E}{\mathbf{E}}
\newcommand{\R}{\mathbb{R}}
\newcommand{\Rp}{\mathbb{R}_+}
\newcommand{\Z}{\mathbb{Z}}
\newcommand{\Zp}{\mathbb{Z}_+}
\newcommand{\Ind}{\mathbf{I}}

\newcommand{\JJ}{J}
\newcommand{\RR}{R}
\renewcommand{\SS}{S}

\renewcommand{\vec}[1]{\text{\boldmath $#1$}}
\newcommand{\modd}[1]{\hat{#1}}
\newcommand{\rs}{{r'}}

\begin{document}

 \begin{center}
   \Large\bfseries Stability of processor sharing networks with
   simultaneous resource requirements
 \end{center}

 \begin{center}
   Jennie Hansen, Cian Reynolds and Stan Zachary
 \end{center}

 \begin{center}
   \textit{Heriot-Watt University\\Edinburgh}
 \end{center}

\begin{quotation}\small\noindent
  We study the phenomenon of \emph{entrainment} in processor sharing
  networks, whereby, while individual network resources have
  sufficient capacity to meet demand, the requirement for simultaneous
  availability of resources means that a network may nevertheless be
  unstable.  We show that instability occurs through poor control, and
  that, for a variety of network topologies, only small modifications
  to controls are required in order to ensure stability.  For controls
  which possess a natural monotonicity property, we give some new
  results for the classification of the corresponding Markov
  processes, which lead to conditions both for stability and for
  instability.
  \vskip 0.3cm
  \noindent
  \emph{Keywords:} processor sharing networks, stability, entrainment.

  \vskip 0.2cm
  \noindent
  \emph{AMS 2000 subject classification:} Primary: 60K20; Secondary: 60K25.

  \vskip 0.2cm
  \noindent
  {\it Short title.\/}  Stability of processor sharing networks

\end{quotation}

\section{Introduction}
\label{sec:introduction}

Modern communications networks, such as the Internet, are able at any
time to share their resources, for example bandwidth, among those
calls or connections currently in progress.  Such calls may require
simultaneous capacity from several resources in the network.  For
example, ``streaming'' applications require, for their duration, a
reservation of bandwidth on each network link over which they connect.
It may then happen that while each resource in the network, considered
in isolation, has sufficient capacity to service the demand placed on
it, the control of the network is such that the requirement for
\emph{simultaneous} availability of capacity ensures that over time
demand cannot be met, and that the network is unstable, that is, that
the number of calls present in the network tends to infinity.  This is
the phenomenon of \emph{entrainment} and has been previously studied
in this context by various authors, notably Bonald and Massouli\'e
(2001), de Veciana \emph{et al} (2001), Kelly and Williams (2004).  In
particular these authors have considered a broad class of so-called
\emph{fair-sharing} control strategies---see below.  Bonald and
Massouli\'e and de Veciana \emph{et al} have shown that here the
instability problems referred to above do not arise, i.e.\ that
provided the various network resources individually have sufficient
capacity, under any fair-sharing control the network will remain
stable.  In the present paper we study the phenomenon of entrainment
in more detail and characterise some of those conditions which give
rise to its occurrence.  In particular we show that it frequently
arises through poor control whenever there are no calls of certain
classes in the network, and that only minor adjustments to control
strategies are required in order to avoid it.  Thus very flexible
management schemes, including those utilising significant
prioritisation, may be safely implemented.  We further establish, in
Section~\ref{sec:stab-prior-based}, some new results for the stability
of Markov chains whose transitions rates possess a natural
monotonicity property, yielding both necessary and sufficient
conditions for the stability of a wide class of network controls.
Since network parameters---in particular call arrival rates---are
unlikely to be known in advance, it is of particular concern to
identify controls whose stability is robust with respect to variations
in these parameters.

We take as our model the following, which is essentially that
introduced by Roberts and Massouli\'e (2000).  Let $\RR$ denote the
finite set
of possible call, or connection, types.  We denote the state of the
network at time~$t$ by $\vec{n}(t)=(n_r(t),\,r\in\RR)$, where $n_r(t)$
is the number of calls of each type~$r$ in progress at that time.
Calls of each type~$r\in\RR$ arrive at the network as a Poisson
process with rate~$\nu_r$ and have \emph{sizes} which are
exponentially distributed with mean~$\mu_r^{-1}$.  Arrival processes
and call sizes are all independent.  (As usual, the above
distributional assumptions make for simplicity of analysis.  However,
it seems likely that all the results of this paper remain
qualitatively correct for other distributions of interarrival times
and call sizes, subject only to these distributions having finite
means, and to the same independence assumptions.  )  When the state of the
network is $\vec{n}$, calls of each type~$r$ are allocated \emph{in
  total} a \emph{bandwidth}~$b_r(\vec{n})\ge0$.  We assume always that
$b_r(\vec{n})=0$ whenever $n_r=0$.  The process~$\vec{n}(\cdot)$ is
thus Markov with state space~$\Zp^{|\RR|}$
and transition
rates given by, for all $r\in\RR$,
\begin{equation}
  \label{eq:1}
  \vec{n} \to
  \begin{cases}
    \vec{n} + \vec{e}_r & \quad\text{at rate $\nu_r$,}\\
    \vec{n} - \vec{e}_r & \quad\text{at rate $\mu_r b_r(\vec{n})$,}
  \end{cases}
\end{equation}
where $\vec{e}_r$ denotes the vector whose $r$th component is $1$
and whose other components are $0$.  For each $\vec{n}$ define also
$\vec{b}(\vec{n})=(b_r(\vec{n}),\,r\in\RR)$.  We shall refer to
$\vec{b}=(\vec{b}(\vec{n}),\,\vec{n}\in\Zp^{|\RR|})$ as the \emph{control}
for the Markov process~$\vec{n}(\cdot)$.  The set of \emph{feasible}
controls~$\vec{b}$ is defined by a set of capacity constraints
\begin{equation}
  \label{eq:2}
  \sum_{r\in\RR}a_{jr}b_r(\vec{n})\le c_j,
  \qquad j\in\JJ,
\end{equation}
indexed in a finite set~$\JJ$, where each $a_{jr}\in\{0,1\}$.  Here
$c_j$ may be thought of as the capacity of resource~$j$, and a
call of type~$r$ utilises resource~$j$ if and only if $a_{jr}=1$.  

In order to allow some results to be stated with sufficient
generality, it is convenient to allow the possibility that, for any
$j$, we may have $c_j=\infty$ (corresponding to the effective
nonexistence of the resource constraint~$j$).  However, we assume,
without loss of generality, that, for all $r\in\RR$, there exists some
$j\in\JJ$ with $c_j<\infty$ and $a_{jr}=1$.

Our interest is in identifying, for fixed values of the
parameters~$\nu_r$, $\mu_r$, $c_j$, and $a_{jr}$, $r\in\RR$,
$j\in\JJ$ (which we shall regard as defining a given network) those
controls~$\vec{b}$ which are both feasible and \emph{stable}, where we
now take the latter to mean that the corresponding Markov
process~$\vec{n}(\cdot)$ is positive recurrent.  A control~$\vec{b}$
for which $\vec{n}(\cdot)$ is null recurrent or transient is referred
to as \emph{unstable}.

For any feasible control~$\vec{b}$, for any $\vec{n}$, we shall say
that a resource~$j\in\JJ$ is \emph{saturated} if the corresponding
constraint~\eqref{eq:2} is satisfied with equality.  Following Bonald
and Massouli\'e (2001), we
shall further say that a feasible control~$\vec{b}$ is \emph{Pareto
  efficient} if, for all $\vec{n}$ and for all $r$ such that $n_r>0$,
there exists $j\in\JJ$ with $a_{jr}=1$ such that $j$ is saturated (so
that $b_r(\vec{n})$ may not be increased without either decreasing
$b_{r'}(\vec{n})$ for some $r'\neq{}r$ or else violating the
constraints~\eqref{eq:2}.)  It is sometimes helpful to consider
controls which are not feasible, so we note in particular that the
requirement of Pareto efficiency includes that of feasibility.

For each $r\in\RR$ we define $\kappa_r=\nu_r/\mu_r$, which may be
thought of as the rate at which ``work'' of type~$r$ arrives at the
network.  Many, but not all, stability results depend on the
parameters $\nu_r$ and $\mu_r$ only through the corresponding
$\kappa_r$.  As observed by Bonald and Massouli\'e (2001), a necessary
and sufficient condition for the existence of a stable feasible
control is given by
\begin{equation}
  \label{eq:3}
  \sum_{r\in\RR}a_{jr}\kappa_r < c_j, \qquad j\in\JJ.
\end{equation}
For the necessity of this condition, observe that if it is violated
for some~$j$, then since, for any feasible control~$\vec{b}$, we then
have $\sum_{r\in\RR}a_{jr}(\kappa_r-b_r(\vec{n}))\ge0$, standard
arguments---see, for example, Asmussen (2003, Chapter~1,
Proposition~5.4) or the argument of Kelly and Williams (2004)---show
that $\vec{b}$ cannot be stable.  For the sufficiency of the
condition~\eqref{eq:3}, define, for any vector
$\hat{\vec{b}}=(\hat{b}_r,\,r\in\RR)$, the corresponding
\emph{complete partitioning} control~$\vec{b}$ by
$b_r(\vec{n})=\hat{b}_r$ whenever $n_r>0$.  Note that this defines a
process~$\vec{n}(\cdot)$ which corresponds to $|\RR|$ independent queues
and is such that, in each queue~$r$, arrivals occur at rate~$\nu_r$
and departures (when $n_r>0$) at rate~$\mu_r\hat{b}_r$.  Thus a
necessary and sufficient condition for the stability of this control
is that $\hat{b}_r>\kappa_r$ for all $r\in\RR$.  In particular we may
choose $\hat{\vec{b}}$ such that the corresponding complete
partitioning control is stable and feasible if and only if the
condition~\eqref{eq:3} holds.  Further, given such a $\hat{\vec{b}}$,
we may clearly define a Pareto efficient control~$\vec{b}'$ such that
$b'_r(\vec{n})\ge\hat{b}_r$ for all $\vec{n}$ and for all $r$ such
that $n_r>0$.  The corresponding process~$\vec{n}'(\cdot)$ may then be
coupled to the process~$\hat{\vec{n}}(\cdot)$ corresponding to
$\hat{\vec{b}}$ in such a way that $n'_r(t)\le\hat{n}_r(t)$ for all
$r$ and for all $t\ge0$.  Hence the condition~\eqref{eq:3} is also
sufficient for the existence of \emph{some} stable Pareto efficient
control.

It is also clear that, in the case $|\JJ|=1$ of a single resource
constraint, the condition~\eqref{eq:3} is sufficient to ensure the
stability of any Pareto efficient control.
This is not in general true when $|\JJ|>1$, as is shown by
Example~\ref{ex:1} below, which is a simplification of one given by
Bonald and Massouli\'e (2001).

\begin{example}\label{ex:1}
  Suppose that $\RR=\{1,2\}$, $\JJ=\{1,2\}$ and that the matrix
  $A=(a_{jr},\,j\in\JJ,r\in\RR)$ is given by
  \begin{displaymath}
    A =
    \begin{pmatrix}
      1 & 1\\
      0 & 1
    \end{pmatrix}
    .
  \end{displaymath}
  Thus calls of type $1$ are constrained by resource $1$ only, while
  calls of type $2$ are constrained by both resources $1$ and $2$.
  Suppose, without loss of generality, 
  that $c_2\le{}c_1$.
  We assume the condition~(\ref{eq:3}), which here becomes
  \begin{equation}
    \label{eq:4}
    \kappa_1+\kappa_2<c_1, \qquad\qquad \kappa_2<c_2.
  \end{equation}
  Consider the Pareto efficient control in which complete priority is
  given to calls of type~$1$.  Then the long-run fraction of time in
  which the network is empty of calls of type~$1$---and so
  resource~$1$ is available for use by calls of type~$2$---is given by
  $1-\kappa_1/c_1$.  Since, when this resource is available, calls of
  type~$2$ are processed at rate~$c_2$, standard arguments for the
  stability of a single-server queue now show that the control is
  stable if and only if $\kappa_2<c_2(1-\kappa_1/c_1)$.  This is a
  condition which is generally more restrictive than \eqref{eq:4}
  above.  When it is violated we have the phenomenon of entrainment
  referred to above; that is, from \eqref{eq:4}, each resource in the
  network, considered in isolation, has sufficient capacity, but the
  given  Pareto efficient control is nevertheless unstable.
\end{example}

In Example~\ref{ex:1} above instability may be considered as occurring
because such high priority is given to calls of type~$1$ as to ensure
that, when these are emptied from the system, resource $1$ is
thereafter consistently underutilised.  As we shall show later, this
problem may in general be avoided by choosing Pareto efficient
controls $\vec{b}$ such that, for each $r$, $b_r(\vec{n})$ is usually
small whenever $n_r$ is small.  In particular this property is
possessed by the class of \emph{(weighted) fair-sharing} controls
introduced by Mo and Walrand (2000) in a generalisation of various
classes considered earlier (see, for example, Kelly \textit{et al}
(1998)).  These are defined as follows: for $\alpha>0$, $\alpha\neq1$,
and weights $w_r$,
$r\in\RR$, a \emph{weighted $\alpha$-fair} control, or bandwidth
allocation, is given by taking, for each $\vec{n}$, $\vec{b}(\vec{n})$
to maximise the concave function
\begin{equation}
  \label{eq:50}
  \sum_{r\in\RR} w_r n_r^\alpha\frac{b_r(\vec{n})^{1-\alpha}}{1-\alpha},
\end{equation}
subject to the capacity constraints~\eqref{eq:2}.  This class is
further extended to each of the cases $\alpha=0,1,\infty$ by taking
the limit of the $\alpha$-fair control as $\alpha$ tends to each of
these values.  (For $\alpha=1$ this is equivalent to replacing the
quantity defined in~\eqref{eq:50} by $\sum_{r\in\RR} w_{r}n_r \log
b_r(\vec{n})$.)  For the cases $\alpha=1$ and $\alpha=\infty$, De
Veciana \textit{et al} (2001)
use Lyapunov function techniques to show that the
condition~\eqref{eq:3} is sufficient for the stability of (weighted)
$\alpha$-fair controls.  Bonald and Massouli\'e (2001) show that this
result holds for general $\alpha$ by using fluid limits and appealing
to a result of Dai (1995) for multiclass queueing networks.  (In fact
it is not certain that Dai's result is directly applicable to the
present networks with their simultaneous resource requirements;
however, Bonald and Massouli\'e's approach is essentially equivalent
to showing that their function~$f$ defined by
\begin{displaymath}
  f(\vec{n})
  = \sum_{r\in\RR} w_r \mu_r^{-1}\kappa_r^{-\alpha}\frac{n_r^{\alpha+1}}{\alpha+1}
\end{displaymath}
is a suitable Lyapunov function for establishing the sufficiency of
\eqref{eq:3} for the stability of $\alpha$-fair controls, and this
requires only a small extension to their argument.)  We note in
particular that the stability of the fair-sharing controls is robust
with respect to parameter variation, subject of course to \eqref{eq:3}
being satisfied.

The resource underutilisation of Example~\ref{ex:1} may be further
understood as resulting from the nonsmooth nature of the call arrival
process.  Consider the analogous \emph{fluid model} in which ``work''
of each type~$r$ arrives steadily at rate $\kappa_r$ and is processed
at rate~$b_r(\vec{n})$, where each $n_r$ is now the volume of work of
type~$r$ in the network and where $\vec{b}$ is again subject to
constraints of the form~(\ref{eq:2}).  Then, under the
condition~\eqref{eq:4} of Example~\ref{ex:1}, it is easy to see that
every Pareto efficient control is stable, in the sense here that the
total volume of work in the system eventually becomes and remains
zero.  For our stochastic model,
the possible modified control discussed above, in which $b_1(\vec{n})$
is kept small whenever $n_1$ is small and $n_2$ is large, may be seen
as a smoothing operation forcing the behaviour of the stochastic model
to follow more closely that of the fluid model.  However, as
Example~\ref{ex:2} below shows, instability may also occur in ways
such that even the corresponding control for the analogous fluid model
is also unstable.

For any control~$\vec{b}$ and any function~$f$ on $\Zp^{|\RR|}$, define
the function $D_\vec{b}f$ on $\Zp^{|\RR|}$ by
\begin{equation}
  \label{eq:5}
  D_\vec{b}f(\vec{n}) = \sum_{r\in\RR}
  \left[\nu_r\left(f(\vec{n}+\vec{e}_r)-f(\vec{n})\right)
    + \mu_rb_r(\vec{n})\left(f(\vec{n}-\vec{e}_r)-f(\vec{n})\right)\right].
\end{equation}
(Since, for $\vec{n}$ and $r$ such that $n_r=0$, we have also
$b_r(\vec{n})=0$, there is no problem arising from the lack of a
formal definition of $f(\vec{n}-\vec{e}_r)$ in this case.  Further,
$D_\vec{b}$ may be thought of as the generator of the Markov
process~$\vec{n}(\cdot)$ under the control~$\vec{b}$.)

\begin{example}\label{ex:2}
  Consider the network defined by $\RR=\{1,2,3\}$,
  $\JJ=\{1,2,3\}$,
  \begin{displaymath}
    A =
    \begin{pmatrix}
      0 & 1 & 1\\
      1 & 0 & 1\\
      1 & 1 & 0
    \end{pmatrix}
    ,
  \end{displaymath}
  with $\nu_r=\nu$, $\mu_r=1$ for all $r$ and $c_j=c$ for all $j$.
  Assume that the condition~\eqref{eq:3} is satisfied, i.e.\ that
  $2\nu<c$.  Consider any control $\vec{b}$ such that, for all
  $\vec{n}\neq\vec{0}$, we have $b_r(\vec{n})=c$ for some $r$ such
  that $n_r>0$ and $b_s(\vec{n})=0$ for $s\ne{}r$.  Then $\vec{b}$ is
  Pareto efficient.  However, for the function
  $f(\vec{n}):=\sum_{r\in\RR}n_r$ and for all
  $\vec{n}\neq\vec{0}$, we have $D_\vec{b}f(\vec{n})=3\nu-c$.
  Standard martingale arguments---again see Asmussen (2003, Chapter~1,
  Proposition~5.4)---now show that $\vec{b}$ is unstable whenever
  $3\nu\ge{}c$.  In this example the instability is not simply the
  result of poor control for $\vec{n}$ close to the boundary of
  $\Zp^3$, and is equally present in the analogous fluid model.  We
  return to this example in Section~\ref{sec:stab-prior-based}.
\end{example}

Our aim in the remainder of the paper is to identify more general classes
of stable controls, to provide insight into how the phenomenon of
entrainment rises, and to show how controls may be modified if
necessary so as to avoid it.  In Section~\ref{sec:workl-based-techn}
we consider a number of fairly simple network topologies, and use
Lyapunov function techniques to show that, for these, any Pareto
efficient control~$\vec{b}$ is stable, provided \emph{only} that, as
discussed above and for suitable $r$, $b_r(\vec{n})$ is modified so as
to be small whenever $n_r$ is small.  In
Section~\ref{sec:stab-prior-based} we consider controls which possess
a natural monotonicity property, likely to be satisfied in any
applications.  We introduce some new analytical techniques to prove
some fairly general results for the stability of the associated Markov
chains.  These give sufficient conditions for stability, which, for
many classes of control, are also close to being necessary.  These
results are applicable to a substantial class of priority-based
controls.

\section{Simple network topologies}
\label{sec:workl-based-techn}

In this section we consider a number of network topologies and in each
case show that, under the condition~\eqref{eq:3}, \emph{any} Pareto
efficient control is stable provided only that it is suitably modified
for values of $\vec{n}$ close to the boundary of $\Zp^{|\RR|}$.

The results of this section are based on simple Lyapunov function
techniques, in particular Foster's criterion.
Proposition~\ref{prop:1} below states the specialisation of this to
the present problem---see, for example, Asmussen (2003, Chapter~1,
Proposition~5.3(ii)), noting that here jumps of the
process~$\vec{n}(\cdot)$ may
only occur between neighbouring states, and also that the usual
uniformisation argument translates statements for discrete-time
processes to the present continuous time setting.

\begin{prop}\label{prop:1}
  Given any control~$\vec{b}$, suppose that there exists a positive
  function~$f$ on $\Zp^{|\RR|}$, a finite subset~$F$ of $\Zp^{|\RR|}$
  and some $\epsilon>0$ such that $D_\vec{b}f(\vec{n})\le-\epsilon$
  for all $\vec{n}\notin F$.
  Then $\vec{b}$ is stable.
\end{prop}

In considering the various network topologies of this section, it is
convenient to define, for any $a\ge1$, the
function~$g_{a}$ on $\Zp$ by
  \begin{displaymath}
    g_{a}(n) =
    \begin{cases}
      \displaystyle{\frac{a}{2} + \frac{n^2}{2a}},
      & \quad\text{if $n < a$}, \\
      n, & \quad\text{if $n \ge a$},
    \end{cases}
  \end{displaymath}  

  We first consider again the very simple Example~\ref{ex:1} of the
  Introduction, but with a general control~$\vec{b}$.

\begin{example}\label{ex:3}
  Let $\RR$, $\JJ$ and $A$ be as in Example~\ref{ex:1} and assume
  again that $c_2\le{}c_1$ and also the condition~\eqref{eq:4}.  We
  show that, given any $\delta>0$, there is some $a\ge1$ (depending on
  the parameters~$\nu_r$, $\mu_r$ and $c_j$ and $\delta$) such that
  a sufficient condition for the stability of any Pareto efficient
  control~$\vec{b}$ is given by
  \begin{equation}
    \label{eq:8}
    \kappa_2 - b_2(\vec{n}) \le -\delta,
    \qquad\text{whenever $n_1<a$ and $\vec{n}\notin{}F$,}
  \end{equation}
  for some finite set~$F$ (necessarily including
  $\{\vec{n}\colon{}n_1<a,n_2=0\}$).  (Note that, from~\eqref{eq:4},
  condition~\eqref{eq:8} may be satisfied for any
  $\delta\in(0,c_2-\kappa_2]$.)

  Observe first that, from \eqref{eq:4} and the saturation of the
  resource~$1$ whenever $\vec{n}\neq\vec{0}$, there exists
  $\delta'>0$ such that, for any Pareto efficient
  control~$\vec{b}$,
  \begin{equation}
    \label{eq:9}
    \sum_{r=1}^2(\kappa_r-b_r(\vec{n})) \le -\delta',
    \qquad\text{for all $\vec{n}\neq\vec{0}$}.
  \end{equation}

  For any fixed $a$, define the function~$f_a$ on $\Zp^2$ by
  \begin{displaymath}
    f_a(\vec{n}) = \frac{g_a(n_1)}{\mu_1} + \frac{n_2}{\mu_2}.
  \end{displaymath}
  Elementary calculations show that, from~\eqref{eq:5}, for any Pareto
  efficient control~$\vec{b}$ satisfying condition~\eqref{eq:8}
  for some finite~$F$, and for any $\vec{n}\notin{}F$,
  \begin{align}
    D_\vec{b}f_a(\vec{n})
    & = \min \left(\frac{n_1}{a},1 \right)[\kappa_1 - b_1(\vec{n})]
    + \kappa_2 - b_2(\vec{n})
    + \frac{1}{2a} h_{1,a}(\vec{n}) \nonumber\\
    & \le - \min \left(\frac{n_1}{a},1 \right)\delta'
    - \max\left(1-\frac{n_1}{a},\,0\right)\delta
    + \frac{1}{2a} h_{1,a}(\vec{n}) \label{eq:10}\\
    & \le - \min(\delta,\delta') + \frac{1}{2a} h_{1,a}(\vec{n}),
    \label{eq:11}
  \end{align}
  where the function $h_{1,a}$ is given by
  \begin{equation}
    \label{eq:12}
    h_{1,a} (\vec{n}) =
    \begin{cases}
      \kappa_1 + b_1(\vec{n}), & \quad\text{if $0 \le n_1 < a$,}\\
      b_1(\vec{n}), & \quad\text{if $n_1 = a$,}\\
      0, & \quad\text{if $n_1 > a$,}
    \end{cases}
  \end{equation}
  and where \eqref{eq:10} follows from \eqref{eq:8} and
  \eqref{eq:9}.  Since $b_1(\vec{n})\le{}c_1$ for any feasible
  control~$\vec{b}$, it follows from \eqref{eq:11} and \eqref{eq:12}
  that $a$ may be chosen sufficiently large that, for any Pareto
  efficient control~$\vec{b}$ satisfying~\eqref{eq:8} for some finite~$F$,
  \begin{displaymath}
    D_\vec{b}f_a(\vec{n}) \le -\frac{1}{2}\min(\delta,\delta')
    \qquad\text{for all $\vec{n}\notin{}F$},
  \end{displaymath}
  and so, by Proposition~\ref{prop:1}, $\vec{b}$ is stable.

  Thus, for this example, any Pareto efficient control is stable
  provided only that, outside of some finite set $F$, it is suitably
  modified for values of $\vec{n}$ such that $n_1$ is small.  However,
  in practice the parameters~$\nu_r$ in particular are unlikely to be
  known, and so it is desirable to choose controls whose stability is
  robust.  One such possibility is to choose any Pareto efficient
  control which assigns complete priority to calls of type~2 whenever
  $n_1<a$ for some $a$.  Provided only that the condition~\eqref{eq:4}
  is satisfied, the above result shows that this will be stable
  provided $a$ is sufficiently large.  The precise value of $a$
  required depends on the slack in the inequalities~\eqref{eq:4}, but
  calculations for the ``worst case'', in which calls of type~1 have
  complete priority whenever $n_1\ge{}a$, show that in general $a$
  need only be small.

  Note in particular that the various fair-sharing controls defined in
  the Introduction always satisfy the condition~\eqref{eq:8} for some
  $\delta$ and for some sufficiently large~$F$, and hence (as already
  remarked) are always stable provided that the condition~\eqref{eq:4}
  is satisfied.
\end{example}

We now extend the above example to each of two more general network
topologies.  In each case it is again the case that only small
modifications, identified below, are required to Pareto efficient
controls in order to ensure their stability.

\begin{example}
  \label{ex:4}
  Consider the network with $\RR=\{1,\dots,|\RR|\}$, and in
  which each call type~$r$ requires service from a dedicated resource
  of capacity $c_r$ together with service from a resource which is
  shared by all call types and has capacity~$c_0$.  The
  constraints~\eqref{eq:2} defining the feasible controls~$\vec{b}$
  are thus
  \begin{align}
    \sum_{r\in\RR} b_r(\vec{n}) & \le c_0, \label{eq:13}\\
    b_r(\vec{n}) & \le c_r, \qquad r\in\RR. \label{eq:14}
  \end{align}

  We assume, without loss of generality, that $c_0<\infty$ and that
  \begin{equation}
    \label{eq:15}
    c_0 \le \sum_{r\in\RR} c_r.
  \end{equation}
  As usual we assume the condition~\eqref{eq:3}, which here becomes
  \begin{align}
    \sum_{r\in\RR} \kappa_r & < c_0, \label{eq:16}\\
    \kappa_r & < c_r, \qquad r\in\RR. \label{eq:17}
  \end{align}
  The case~$|\RR|=2$ with $c_1=\infty$ is the earlier
  Example~\ref{ex:3} (with $c_0$ here corresponding to $c_1$ there).
  However, in the general case described above a little more care is
  needed in the conditions for the stability of Pareto efficient
  controls.

  For any $\vec{n}\in\Zp^{|\RR|}$, define $n_{\min}=\min_{r\in\RR}n_r$.  We
  show briefly that, given any $\delta>0$, there is again some $a\ge1$
  (depending on the parameters~$\nu_r$, $\mu_r$ and $c_j$ and
  $\delta$) such that a sufficient condition for the stability of any
  Pareto efficient control~$\vec{b}$ is given by
  \begin{equation}
    \label{eq:18}
      \sum_{r\colon n_r\ge a'}  (\kappa_r - b_r(\vec{n})) \le -\delta,
    \qquad\text{for all $a'\in[1,a]$ and $\vec{n}$
      such that $n_{\min}<a$, $\vec{n}\notin{}F$,}
  \end{equation}
  for some finite set~$F$.  (The existence of a $\delta>0$ such that
  the condition~\eqref{eq:18} may be satisfied is guaranteed by
  \eqref{eq:16} and \eqref{eq:17}.  Further, it is not difficult to
  see that in the case $R=2$, this condition reduces to that of
  Example~\ref{ex:3} taken together with a similar condition in which
  the roles of calls of types~$1$ and $2$ are interchanged.)

  The proof of this result is similar to that of Example~\ref{ex:3}.
  Note first that it follows from~(\ref{eq:15}) that, for any Pareto
  efficient control, the resource~$0$ is necessarily saturated for any
  $\vec{n}$ such that $n_{\min}\ge1$.  Hence, from~(\ref{eq:17}),
  there exists $\delta'>0$ such that, again for any Pareto
  efficient control~$\vec{b}$,
  \begin{equation}
    \label{eq:19}
    \sum_{r\in\RR}(\kappa_r - b_r(\vec{n})) \le -\delta'
    \qquad\text{for all $\vec{n}$ such that $n_{\min}\ge1$.}
  \end{equation}
  For any fixed $a$, define the function~$f_a$ on $\Zp^{|\RR|}$ by
  \begin{displaymath}
    f_a(\vec{n}) = \sum_{r\in\R}\frac{g_a(n_r)}{\mu_r}.
  \end{displaymath}
  Elementary calculations show that, from~\eqref{eq:5}, for any Pareto
  efficient control~$\vec{b}$ satisfying the condition~\eqref{eq:18}
  for some finite~$F$, and for any $\vec{n}\notin{}F$,
  \begin{align}
    D_\vec{b}f_a(\vec{n})
    & = \sum_{r\in\RR}%
    \min \left(\frac{n_r}{a},1 \right)[\kappa_r - b_r (\vec{n})]
    + \frac{1}{2a} \sum_{r\in\RR} h_{r,a}(n_r),
    \nonumber\\
    & = \frac{1}{a}\sum_{a'=1}^{a}\;\sum_{r\colon n_r\ge a'}%
    (\kappa_r - b_r(\vec{n}))
    + \frac{1}{2a} \sum_{r\in\RR} h_{r,a}(n_r) \nonumber\\
    & \le -\min(\delta,\delta') + \frac{1}{2a} \sum_{r\in\RR} h_{r,a}(n_r),
    \label{eq:20}
  \end{align}
  where, for each $r$, the function~$h_{r,a}$ is given
  by~(\ref{eq:12}) with the index $r$ replacing the index~$1$, and
  where \eqref{eq:20} follows from \eqref{eq:18} for $\vec{n}$ such
  that $n_{\min}<a$ and from \eqref{eq:19} for $\vec{n}$ such that
  $n_{\min}\ge{}a$.  It now follows as in Example~\ref{ex:3} that $a$
  may be chosen sufficiently large that, for any Pareto efficient
  control~$\vec{b}$ satisfying~\eqref{eq:18} for some finite~$F$,
  \begin{displaymath}
    D_\vec{b}f_a(\vec{n}) \le -\frac{1}{2}\min(\delta,\delta')
    \qquad\text{for all $\vec{n}\notin{}F$},
  \end{displaymath}
  and so, again by Proposition~\ref{prop:1}, the control~$\vec{b}$ is
  stable.

  Thus, in order to ensure the stability of a general Pareto efficient
  control~$\vec{b}$, it is only necessary to appropriately modify
  $\vec{b}(\vec{n})$ for $\vec{n}$ such that $n_{\min}<a$.  In
  particular it follows from the above result that a Pareto efficient
  control~$\vec{b}$ whose stability is reasonably robust is given by
  requiring that, for some~$a$, for all $\vec{n}$ (such that
  $n_{\min}<a$) and for all $a'\in[1,a]$, calls of types~$r$ such that
  $n_r\ge{}a'$ collectively have complete priority over calls of the
  remaining types, that is, that
  \begin{displaymath}
    \sum_{r\colon n_r\ge a'} b_r(\vec{n})
    = \min\biggl(c_0, \sum_{r\colon n_r\ge a'} c_r\biggr).
  \end{displaymath}
  For $a$ sufficiently large, depending on the slack in the
  inequalities~(\ref{eq:16}) and (\ref{eq:17}), any such control is
  stable.

  Now consider further the case~$R=2$.  Note that the present topology
  is completely general for a network with two call types.
  Let~$\vec{b}$ be any Pareto efficient control such that
  \begin{align}
    \lim_{n_1\to\infty}b_1(n_1,n_2) & = c_1
    \qquad\text{for all $n_2$}, \label{eq:21}\\
    \lim_{n_2\to\infty}b_2(n_1,n_2) & = c_2
    \qquad\text{for all $n_1$}. \label{eq:22}
  \end{align}
  Then it is straightforward that, for
  $\delta<\min_{r=1,2}(c_r-\kappa_r)$, and for any $a$, the
  condition~\eqref{eq:18} is satisfied for $F$ sufficiently large, and
  so $\vec{b}$ is stable.

  In particular the conditions~\eqref{eq:21} and \eqref{eq:22} are
  again satisfied by the various fair-sharing controls.  The present
  conditions are of course considerably more general.  However the
  development of corresponding results for general networks with
  $R\ge3$ remains a challenging problem.
\end{example}

Our final network topology is a simple ``backbone'' structure.

\begin{example}\label{ex:5}
  Consider a network with resource set $\JJ=\{1,\dots,k\}$, where
  resource~$j$ has capacity~$c_j$ as usual.  The set of call types is
  given by $\RR=\{0,1,\dots,k\}$ where calls of each
  type~$r=1,\dots,k$ require service from the single resource~$j=r$,
  while calls of type~$0$ require service from each of the
  resources~$1,\dots,k$.

  The constraints~\eqref{eq:2} defining the feasible
  controls~$\vec{b}$ are given by
  \begin{equation}\label{eq:23}
    b_0(\vec{n}) + b_r(\vec{n}) \le c_r,
    \qquad r=1,\dots k.
  \end{equation}
  Again we assume the condition~\eqref{eq:3}, which here becomes
  \begin{equation}\label{eq:24}
    \kappa_0 + \kappa_r  < c_r,
    \qquad r=1,\dots k.
  \end{equation}
  The state of the network is thus here denoted by
  $\vec{n}=(n_0,n_1,\dots,n_k)$.  For any such $\vec{n}$ define
  $\hat{n}_{\max}=\max(n_1,\dots,n_k)$.  
  
  This example again generalises that of Example~\ref{ex:3} with the
  call type~$0$ here playing the role of the call type~$2$ in that
  example.  In a fairly straightforward generalisation of the result
  of that example, it is here the case that, given any $\delta>0$,
  there is again an $a\ge1$ such that a sufficient condition for the
  stability of any Pareto efficient control~$\vec{b}$ is given by
  \begin{equation}
    \label{eq:55}
    \kappa_0 - b_0(\vec{n}) \le -\delta,
    \qquad\text{whenever $\hat{n}_{\max}<a$ and $\vec{n}\notin{}F$,}
  \end{equation}
  for some finite set~$F$

  To see this, suppose first that $\mu_r=1$ for all $r=0,1,\dots,k$
  and that the condition~\eqref{eq:55} is satisfied.
  For any $a\ge1$, define the function~$f_a$ on $\Zp^{|\RR|}$ by 
  \begin{math}
    f_a(\vec{n}) = n_0 + g_a(\hat{n}_{\max}).
  \end{math}
  Then, as in the case of Example~\ref{ex:3}, it follows
  straightforwardly from the conditions~\eqref{eq:24} and
  \eqref{eq:55} that there exists $\delta'>0$ and $a$ sufficiently
  large that $D_\vec{b}f_a(\vec{n})\le-\delta'$ for all
  $\vec{n}\in\Zp^{|\RR|}\setminus{}F$ such that additionally
  $\hat{n}_{\max}=n_r$ for a \emph{single} value of $r=1,\dots,k$.  It is now
  easy to see that this latter restriction may be removed (possibly at
  the expense of increasing $a$) by suitably smoothing the
  function~$f$ in the neighbourhood of those $\vec{n}$ such that
  $\hat{n}_{\max}=n_r$ for two or more values of $r=1,\dots,k$.  The
  desired result thus follows in this case; for general $\mu_r$ only
  routine modifications to the above argument are required.

  Hence we again have that any Pareto efficient control requires only
  slight modification---for those $\vec{n}$ such that
  $\hat{n}_{\max}<a$---in order to be stable.  A robust Pareto
  efficient control is given, for example, by assigning complete
  priority to calls of type~$0$ whenever $\hat{n}_{\max}<a$, the
  necessary value of $a$ depending on the slack in the
  inequalities~\eqref{eq:24}.
\end{example}

Examples~\ref{ex:3}--\ref{ex:5} above all consider fairly simple
network topologies.  In the analogous fluid model defined in the
Introduction, it is easily seen that, for each of these topologies,
condition~\eqref{eq:3} is sufficient for the stability of \emph{any}
Pareto efficient control.  (In each case this follows, for example, by
using the same Lyapunov function as for the stochastic model, except
that the function $g_a$ may be replaced by the identity function.)
The examples illustrate a principle which seems likely to be true for
more general network topologies, namely that when a control is such
that it is stable for the fluid model, then there is a closely
approximating control which is stable for the corresponding stochastic
model.

In the next section we consider stability criteria for quite general
network topologies, applicable typically to controls where there is
some prioritisation among call types.

\section{Stability of monotonic controls}
\label{sec:stab-prior-based}

Many controls likely to be of practical application possess a
simple monotonicity property (see below).  In this section we study
stability for a wide class of such controls, giving sufficient
conditions for stability, which, for many classes of control, are also
close to being necessary (see Remark~\ref{rem:stab} below which
further discusses the applicability of the results of this section).
We require first the following quite general lemma.

We shall say that a control~$\vec{b}$ is \emph{bounded} if, for all
$r$, $b_r(\vec{n})$ is bounded in $\vec{n}$.

\begin{lem}\label{lem:over}
  Let $\vec{b}$ be any bounded control and, as usual, let
  $\vec{n}(\cdot)$ denote the Markov process for the corresponding
  network.  Then, for all $r\in\RR$,
  \begin{equation}\label{eq:27}
    \limsup_{t\to\infty} \frac{1}{t} \int_{0}^{t}b_r(\vec{n}(u))\,du
    \le \kappa_r \quad\text{a.s.}
 \end{equation}
\end{lem}

\begin{proof}
  For all $r\in\RR$, the compensated process~$n^*_r(\cdot)$ defined by
  \begin{displaymath}
    n^*_r(t) =  n_r(t) - n_r(0) + \int_{0}^{t} (b_r(\vec{n}(u)) - \kappa_r)\, du
  \end{displaymath}
  is a zero-mean martingale.  Since also the transition rates of the
  Markov process~$\vec{n}(\cdot)$ are bounded, it follows that, for
  some constant~$M$ and for all $t\ge0$, we have
  $\E(n^*_r(t)^2)\le{}Mt$.  Thus, for $1/2<\alpha<1$, the process
  \begin{math}
    \left(n^*_r(t)/t^{\alpha} \right)_{t>0}
  \end{math}
  is an $\mathcal{L}^2$-bounded supermartingale, and so, as
  $t\to\infty$, converges almost surely to some finite random
  variable.  It follows that
  \begin{displaymath}
    \frac{n_r(t)}{t} + \frac{1}{t} \int_{0}^{t} b_r(\vec{n}(u))\,du
    \to \kappa_r \quad\text{a.s.}
  \end{displaymath}
  Since the process~$n_r(\cdot)$ is positive, the result~\eqref{eq:27}
  now follows.
\end{proof}

\begin{rem}\label{rem:pr}
  In the case where the bounded control~$\vec{b}$ is stable, the
  Markov process~$\vec{n}(\cdot)$ is positive recurrent, and so we
  have the stronger result that, for all $r\in\RR$,
  \begin{equation}
    \label{eq:28}
    \lim_{t\to\infty} \frac{1}{t} \int_{0}^{t}b_r(\vec{n}(u))\,du
    = \E_\pi b_r = \kappa_r \quad\text{a.s.}
  \end{equation}
  Here $\E_\pi b_r$ denotes the expectation of the function~$b_r$ with
  respect to the stationary distribution~$\pi$ of $\vec{n}(\cdot)$,
  the first equality in \eqref{eq:28} follows from the ergodic
  theorem, and the second is simply the assertion that, under
  stationarity (and readily deducible from the balance equations
  defining~$\pi$), the expected arrival and departure rates for calls
  of type~$r$ are equal.
\end{rem}

We shall say that a bounded control~$\vec{b}$ is \emph{monotonic} if, for all
$r$ and for all $\vec{n}$,
\begin{gather}
  \text{$b_r(\vec{n})$ is increasing in $n_r$ (with $n_s$ fixed for
    all $s\ne{}r$)},  \label{eq:29}\\
  \text{$b_r(\vec{n})$ is decreasing in $n_s$ (with
    $n_{s'}$ fixed for all $s'\ne{}s$), for all $s\ne{}r$.} \label{eq:30}
\end{gather}
Note that, depending on the network structure, this property is
natural in many applications.  For instance, for the structure of
Example~\ref{ex:4}, it is possessed by all the fair-sharing controls,
and also by any other reasonable Pareto efficient control.  For more
complex network structures, controls may be coupled to monotonic
controls to establish stability results using the results given below.
See also Example~\ref{ex:6} below.  We note further that a related but
somewhat different definition of monotonicity is used by Bonald and
Prouti\`{e}re (2004).

For any monotonic control~$\vec{b}$, and for each $\SS\subseteq{}\RR$
(including the case where $\SS$ is the empty set~$\emptyset$), define
the function $\vec{b}^\SS\colon{}\Zp^{|\SS|}\to\Rp^{|\RR|}$ by
\begin{align}
  b^\SS_r(\vec{n}_\SS)
  & = \lim_{n_s\to\infty\,\forall\, s\notin \SS} b_r(\vec{n}),
  \qquad r\in \SS,
  \label{eq:31}\\
  b^\SS_r(\vec{n}_\SS)
  & = \lim_{n_r\to\infty}\,%
  \lim_{n_s\to\infty\,\forall\, s\notin \SS\cup\{r\}} b_r(\vec{n}),
  \qquad r\notin \SS,
  \label{eq:32}
\end{align}
where $\vec{n}_\SS=(n_s,\,s\in{}\SS)$ and where, in taking the limits
in the right side of each of \eqref{eq:31} and \eqref{eq:32}, the
vector of those coordinates of $\vec{n}$ that belong to $\SS$ is held
fixed at $\vec{n}_\SS$.  Note that, by monotonicity, this
function is well-defined.  In particular, in \eqref{eq:31}, the order
within $\RR\setminus{}\SS$ in which the limits are taken is irrelevant;
however, in \eqref{eq:32}, the final limit to be taken must be that as
$n_r\to\infty$.

For any monotonic control~$\vec{b}$ and any $\SS$ as above, we shall
say that $\vec{b}^\SS$ is \emph{stable} if the application of the
control~$(b^\SS_s,\,s\in{}\SS)$ to calls in the network of
types~$s\in{}\SS$ (with $\nu_s$, $\mu_s$, $s\in{}\SS$, as usual)
yields a positive recurrent Markov
process~$\vec{n}^\SS(\cdot)=(n^\SS_s(\cdot),\,s\in{}\SS)$.  (This
$|\SS|$-dimensional process may be thought of as that which results
when the number of calls of each type $r\notin{}\SS$ is infinite.)
When $\vec{b}^\SS$ is stable we shall let $\pi^\SS$ denote the
stationary distribution of $\vec{n}^\SS(\cdot)$ (or, where necessary,
the probability function of this distribution); we shall further
define, for each $r\in\RR$,
\begin{equation}\label{eq:33}
  \E_{\pi^\SS}b^\SS_r
  = \sum_{\vec{n}_\SS\in\Zp^{|\SS|}}\pi^\SS(\vec{n}_\SS) b^\SS_r(\vec{n}_\SS)
\end{equation}
to be the expected value of $b^\SS_r$ under this distribution.  In the
case where $\SS$ is the empty set~$\emptyset$, we have that
$\vec{b}^\emptyset=(b^\emptyset_r,\,r\in\RR)$ is a vector of
constants.  We make the natural convention that $\vec{b}^\emptyset$ is
always stable; the distribution $\pi^\emptyset$ is concentrated on a
single point, and we have
$\E_{\pi^\emptyset}b^\emptyset_r=b^\emptyset_r$ for all $r\in\RR$.

Our main result of this section is Theorem~\ref{thm:priority} below.
The first part is similar in spirit to results of Borovkov (1998,
Chapter~8)
for asymptotically spatially homogeneous Markov processes.  However,
the application of those results here would require that the right
side of~\eqref{eq:32} is invariant under interchange of the limits in
that expression, a condition which is not in general satisfied for our
monotonic controls.  Rather the monotonicity itself provides
sufficient structure to obtain the results of
Theorem~\ref{thm:priority}. 

\begin{theorem}\label{thm:priority}
  Suppose that the control~$\vec{b}$ is monotonic and that
  $\SS\subseteq{}\RR$ is such that $\vec{b}^\SS$ is stable.
  \begin{itemize}\itemsep 0pt
  \item[(i)] If $r\notin{}\SS$ is such that
    \begin{equation}
      \label{eq:34}
      \E_{\pi^\SS}b^\SS_r>\kappa_r,
    \end{equation}
    then $\vec{b}^{\SS\cup\{r\}}$ is stable. 
  \item[(ii)] If $r\notin{}\SS$ is such that
    \begin{equation}
      \label{eq:35}
      \E_{\pi^\SS}b^\SS_r<\kappa_r,
    \end{equation}
    then $\vec{b}^{\SS\cup\{r\}}$ is unstable. 
  \end{itemize}
\end{theorem}

\begin{rem}\label{rem:stab}
  Given the stability of $\vec{b}^\SS$ for some $\SS\subset\RR$
  (recall that, as already remarked, $\vec{b}^\emptyset$ is always
  stable), Theorem~\ref{thm:priority} gives criteria for determining
  the stability or otherwise of $\vec{b}^{\SS\cup\{r\}}$ for any
  $r\notin{}\SS$, except only in the case of equality in~\eqref{eq:34}
  or \eqref{eq:35} (where the natural conjecture is that
  $\vec{b}^{\SS\cup\{r\}}$ is unstable---see also the remarks at the
  end of Example~\ref{ex:6}).  Recursive application of the theorem
  thus yields sufficient conditions both for the stability and the
  instability of monotonic controls.  However, note that, for example
  in the case $\RR=\{1,2\}$, $\vec{b}^{\{1\}}$ and $\vec{b}^{\{2\}}$
  may both be unstable, while $\vec{b}=\vec{b}^{\{1,2\}}$ is stable,
  as is the case for fair-sharing controls here.  In such
  circumstances Theorem~\ref{thm:priority} does not settle the
  question of the stability of $\vec{b}$.  Rather its primary
  application is to controls in which there is a sufficient hierarchy
  of prioritisation as to permit the recursive application of the
  first part of the theorem to at least establish the stability of
  $\vec{b}^\SS$ for $\SS=\RR\setminus\{r\}$ for some $r\in\RR$.  The
  theorem then also (in general) settles the question of the stability
  of $\vec{b}$ itself.  For an illustration of the application of the
  theorem, see Example~\ref{ex:6} below.
\end{rem}

\begin{proof}[Proof of Theorem~\ref{thm:priority}]
  Since, for given $\SS\subset\RR$ and $r\notin{}\SS$, the stability of
  $\vec{b}^{\SS\cup\{r\}}$ corresponds to the positive recurrence of the
  Markov process~$\vec{n}^{\SS\cup\{r\}}(\cdot)$ defined above (in which
  the number of calls of each type~$s\notin{}\SS\cup\{r\}$ is
  effectively held at infinity), it is sufficient to prove the
  results~(i) and (ii) in the case where $\SS=\RR\setminus\{\rs\}$ for
  some~$\rs$, and we henceforth assume this.  (The primary advantage of
  doing so is that we avoid some unpleasant notational complexity.)
  We identify any $\vec{n}\in\Zp^{|\RR|}$ with the pair~$(\vec{n}_\SS,n_\rs)$
  where $\vec{n}_\SS=(n_s,\,s\in{}\SS)$.  Recall that then, for each such
  $\vec{n}_\SS$, we have
  $b^\SS_\rs(\vec{n}_\SS)=\lim_{n_\rs\to\infty}b_\rs(\vec{n}_\SS,n_\rs)$.

  Suppose first that the condition~(\ref{eq:34}) holds.  We require to
  show that $\vec{b}$ is stable.  The underlying idea is that the
  monotonicity of $\vec{b}$ and stability of $\vec{b}^\SS$ ensure that
  the components $(n_s(\cdot),\,s\in{}\SS)$ of the
  process~$\vec{n}(\cdot)$ become and remain small, and the
  condition~\eqref{eq:34} then ensures that, except in some finite
  region~$A$, the remaining component~$n_\rs(\cdot)$ of this process is
  decreasing at a rate which is bounded away from zero; thus the
  process~$\vec{n}(\cdot)$ spends, in the long term, a nonzero
  proportion of time within~$A$.

  It follows from~(\ref{eq:34}) and the monotonicity of $\vec{b}$ that
  we can choose a finite set
  $A=\{\vec{n}\colon{}n_r\le{}\bar{n}_r,\,r\in\RR\}\subset\Zp^{|\RR|}$ and a
  positive function~$\bar{b}_\rs$ on $\Zp^{|\SS|}$ such that
  \begin{align}
    \text{$\bar{b}_\rs$ is decreasing in each of its arguments},
    \hspace{-6em}\label{eq:36}\\
      \sum_{\vec{n}_\SS\in\Zp^{R-1}}\pi^\SS(\vec{n}_\SS)\bar{b}_\rs(\vec{n}_\SS)
      & > \kappa_\rs,  \label{eq:37}\\
      b_\rs(\vec{n}) & \ge \bar{b}_\rs(\vec{n}_\SS)
      \qquad\text{for all $\vec{n}\notin A$}.\label{eq:38}
  \end{align}
  For example, given $\bar{n}_r$, $r\in\RR$, we may choose
  \begin{equation}\label{eq:39}
    \bar{b}_\rs(\vec{n}_\SS)
    = \Ind(n_s\le\bar{n}_s,\,s\in{}\SS)\,b_\rs(\vec{n}_\SS,\bar{n}_\rs),
    \qquad \vec{n}_\SS \in \Zp^{|\SS|},
  \end{equation}
  where $\Ind$ denotes the indicator function.
  The condition~(\ref{eq:36}) then follows from the monotonicity of
  $\vec{b}$; since also $\vec{b}$ is nonnegative the
  condition~(\ref{eq:38}) follows trivially from~(\ref{eq:39}), except
  for $\vec{n}$ such that $n_s\le\bar{n}_s$ for $s\in{}\SS$ and
  $n_\rs>\bar{n}_\rs$, in which case (\ref{eq:38}) again follows from the
  monotonicity of $\vec{b}$; finally the condition~(\ref{eq:37})
  follows from~\eqref{eq:34} and the monotone convergence theorem by
  choosing $\bar{n}_r$, $r\in\RR$, all sufficiently large, since, for
  any $\vec{n}_\SS$,
  \begin{displaymath}
    \lim_{\bar{n}_r\to\infty\,\forall\,r\in\RR}%
    \Ind(n_s\le\bar{n}_s,\,s\in{}\SS)\,b_\rs(\vec{n}_\SS,\bar{n}_\rs)
    = \lim_{\bar{n}_\rs\to\infty}b_\rs(\vec{n}_\SS,\bar{n}_\rs)
    = b^\SS_\rs(\vec{n}_\SS).
  \end{displaymath}
  It also follows from the monotonicity of the control~$\vec{b}$ that
  we can couple the corresponding Markov process~$\vec{n}(\cdot)$ to a
  process~$\vec{n}^\SS(\cdot)=(n^\SS_s(\cdot),\,s\in{}\SS)$ on $\Zp^{|\SS|}$
  with control $\vec{b}^\SS$ in such a way that $n_s(t)\le{}n^\SS_s(t)$
  for all $t>0$ and for all $s\in{}\SS$.  Since also the
  process~$\vec{n}^\SS(\cdot)$ has stationary distribution~$\pi^\SS$, it
  follows from~(\ref{eq:36}) that, for the function~$\bar{b}_\rs$ defined
  above,
  \begin{equation}
    \label{eq:40}
    \limsup_{t\to\infty}\frac{1}{t}\int_0^t\bar{b}_\rs(\vec{n}_\SS(u))\,du
    \ge 
    \lim_{t\to\infty}\frac{1}{t}\int_0^t\bar{b}_\rs(\vec{n}^\SS(u))\,du
    = \E_{\pi^\SS}\bar{b}_\rs,
  \end{equation}
  where $\vec{n}_\SS(\cdot)=(n_s(\cdot),\,s\in{}\SS)$ and where the
  final equality above follows by the ergodic theorem.

  It now follows from Lemma~\ref{lem:over} that, for some $M>0$,
  \begin{align}
    \kappa_\rs
    & \ge \limsup_{t\to\infty}\frac{1}{t}\int_0^t
    b_\rs(\vec{n}(u))\,du
    \nonumber \\
    & \ge \limsup_{t\to\infty}\frac{1}{t}\int_0^t
    b_\rs(\vec{n}(u))\,\Ind((\vec{n}(u)\notin A)\,du
    \nonumber \\
    & \ge \limsup_{t\to\infty}\frac{1}{t}\int_0^t
    \bar{b}_\rs(\vec{n}_\SS(u))\,\Ind((\vec{n}(u)\notin A)\,du
    \label{eq:41}\\
     & \ge \limsup_{t\to\infty}\frac{1}{t}\int_0^t
     \bar{b}_\rs(\vec{n}_\SS(u))\,du
     - \lim_{t\to\infty}\frac{M}{t}\int_0^t\Ind((\vec{n}(u)\in{}A)\,du
     \label{eq:42}\\
     & \ge \E_{\pi^\SS}\bar{b}_\rs
     - \lim_{t\to\infty}\frac{M}{t}\int_0^t\Ind((\vec{n}(u)\in{}A)\,du,
     \label{eq:43}
  \end{align}
  where the inequality~\eqref{eq:41} follows from~\eqref{eq:38}, the
  inequality~(\ref{eq:42}) follows since, from~~\eqref{eq:36}, the
  function~$\bar{b}_\rs$ is necessarily bounded (note that, since $A$ is
  finite, the limit in the final term in \eqref{eq:42} always exists),
  and the final inequality~(\ref{eq:43}) follows from~(\ref{eq:40}).
  Thus, from \eqref{eq:43} and~(\ref{eq:37}),
  \begin{displaymath}
    \lim_{t\to\infty}\frac{1}{t}\int_0^t\Ind((\vec{n}(u)\in{}A)\,du > 0.
  \end{displaymath}
  Since $A$ is finite it now follows from the ergodic theorem that
  the Markov process~$\vec{n}(\cdot)$ is positive recurrent and so
  $\vec{b}$ is stable.

  Now suppose instead that the condition~\eqref{eq:35} holds.  We show
  that the Markov process~$\vec{n}(\cdot)$ corresponding to $\vec{b}$
  is transient (and hence $\vec{b}$ is unstable).  The underlying idea
  here is that whenever $n_\rs(\cdot)$ is very large, the process
  $\vec{n}(\cdot)$ again behaves approximately as if it were
  controlled by $\vec{b}^\SS$, and thus, from~\eqref{eq:35}, we may
  expect that $\lim_{t\to\infty}n_\rs(t)=\infty$ a.s.  To make this
  rigorous we again use the monotonicity of $\vec{b}$ to couple the
  process~$\vec{n}(\cdot)$ to a process~$\modd{\vec{n}}(\cdot)$ whose
  control is sufficiently close to that of $\vec{b}^\SS$ that we may
  show that $\lim_{t\to\infty}\modd{n}_\rs(t)=\infty$ a.s., and for
  which the coupling ensures that also
  $\lim_{t\to\infty}n_\rs(t)=\infty$ with strictly positive probability.

  Given $\bar{n}_\rs\in\Zp$ (fixed, to be chosen later), define a Markov
  process~$\modd{\vec{n}}(\cdot)=(\modd{n}_r(\cdot),\,r\in\RR)$ as
  follows: for each $s\in{}\SS$, the component
  process~$\modd{n}_s(\cdot)$ has state space~$\Zp$ as usual, while
  $\modd{n}_\rs(\cdot)$ has state space~$\Z$; for each $r\in\RR$,
  transitions~$n_r\to{}n_r+1$ occur at rate~$\nu_r$ as usual, while
  transitions~$n_r\to{}n_r-1$ occur at
  rate~$\mu_r\modd{b}_r(\vec{n})$, where
  \begin{align}
    \modd{b}_s(\vec{n}) & = b_s(\vec{n}_\SS,\bar{n}_\rs),
    \qquad s\in \SS, \label{eq:48}\\
    \modd{b}_\rs(\vec{n}) & = b^\SS_\rs(\vec{n}_\SS). \label{eq:49}
  \end{align}
  Observe that the process~$\modd{\vec{n}}(\cdot)$ has uniformly
  bounded transition rates which are independent of $n_\rs\in\Z$.
  Suppose that 
  \begin{equation}
    \label{eq:57}
    \vec{n}(0)=\modd{\vec{n}}(0),  \qquad n_\rs(0)>\bar{n}_\rs.
  \end{equation}
  Define the random time~$T=\min\{t>0\colon{}n_\rs(t)<\bar{n}_\rs\}$.  It
  follows from the monotonicity of $\vec{b}$ that, for all $\vec{n}$
  such that $n_\rs\ge\bar{n}_\rs$,
  \begin{equation}
    \label{eq:45}
    \modd{b}_s(\vec{n})\ge{}b_s(\vec{n})\ge{}b^\SS_s(\vec{n}_\SS),
    \qquad s\in \SS.
  \end{equation}
  and hence that we may couple the processes~$\modd{\vec{n}}(\cdot)$,
  $\vec{n}(\cdot)$ and the $|\SS|$-dimensional
  process~$\vec{n}^\SS(\cdot)$ with control~$(b^\SS_s,\,s\in{}\SS)$ in
  such a way that
  \begin{equation}
    \label{eq:46}
      \modd{n}_s(t) \le n^\SS_s(t),
      \qquad\text{for all $s\in \SS$ and for all $t\ge0$}
    \end{equation}
  and
  \begin{equation}
    \label{eq:25}
    \modd{n}_s(t) \le n_s(t),
    \qquad\text{for all $s\in \SS$ and for all $0\le t\le T$.}
  \end{equation}
  Since $\vec{n}^\SS(\cdot)$ is assumed positive recurrent with
  stationary distribution~$\pi^\SS$, it follows from \eqref{eq:45} or
  \eqref{eq:46} that
  $\modd{\vec{n}}_\SS(\cdot)=(\modd{n}_s(\cdot),\,s\in\SS)$ is
  similarly positive recurrent with stationary
  distribution~$\modd{\pi}_\SS$ say.  Further, as
  $\bar{n}_\rs\to\infty$, the control~$\modd{\vec{b}}$ converges
  pointwise in each of its components to $\vec{b}^\SS$.  Hence
  elementary arguments (e.g.\ consideration of the times of return to
  $\vec{0}$ of the process~$\vec{n}^\SS(\cdot)$, coupled with the use
  of~\eqref{eq:46} and the ergodic theorem) show that, again as
  $\bar{n}_\rs\to\infty$, $\modd{\pi}_\SS$ converges in distribution
  to $\pi^\SS$.  Since also $b^\SS_\rs$ is bounded, it now follows
  from \eqref{eq:35} that we may choose the constant~$\bar{n}_\rs$
  sufficiently large that
  \begin{equation}
    \label{eq:47}
          \E_{\modd{\pi}_\SS}b^\SS_\rs<\kappa_\rs,
  \end{equation}
  (where, analogously to \eqref{eq:33}, $\E_{\modd{\pi}_\SS}b^\SS_\rs$ is the
  expectation of $b^\SS_\rs$ with respect to $\modd{\pi}_\SS$).

  It further follows from the monotonicity of $\vec{b}$ and
  from~\eqref{eq:25} that, for all $0\le t\le T$,
  \begin{align}
    \modd{b}_\rs(\modd{\vec{n}}(t))
    & = b^\SS_\rs(\modd{\vec{n}}_\SS(t)) \nonumber\\
    & \ge b^\SS_\rs(\vec{n}_\SS(t)) \label{eq:26}\\
    & \ge b_\rs(\vec{n}(t)) \label{eq:44},
  \end{align}
  where \eqref{eq:26} follows from~\eqref{eq:25} and the monotonicity
  of $\vec{b}$ while \eqref{eq:44} follows from the definition of
  $\vec{b}^\SS$ and, again, the monotonicity of $\vec{b}$.  It now
  follows from \eqref{eq:44} that we may couple also the
  components~$\rs$ of the processes~$\modd{\vec{n}}(\cdot)$ and 
  $\vec{n}(\cdot)$ in such a way that
  \begin{equation}\label{eq:56}
    \modd{n}_\rs(t) \le n_\rs(t),
    \qquad\text{for all $0\le t\le T$.}
  \end{equation}
  As noted above, the
  process~$\modd{\vec{n}}_\SS(\cdot)=(\modd{n}_s(\cdot),\,s\in{}\SS)$ has
  stationary distribution~$\modd{\pi}_\SS$, while the
  process~$\modd{n}_\rs(\cdot)$ may be viewed as a Markov additive
  process modulated by the remaining components
  $\modd{\vec{n}}_\SS(\cdot)$ of $\modd{\vec{n}}(\cdot)$.
  From~\eqref{eq:47}, the expectation of the increments of
  $\modd{n}_\rs(\cdot)$ between those times at which
  $\modd{\vec{n}}(\cdot)$ returns to any fixed state is strictly
  positive.  It follows from the standard theory of Markov additive
  processes that $\lim_{t\to\infty}\modd{n}_\rs(t)=\infty$ a.s., and
  further that, under the condition~\eqref{eq:57},
    \begin{displaymath}
      \Pr(\modd{n}_\rs(t)\ge\bar{n}_\rs \text{ for all $t\ge0$,}\ %
      \lim_{t\to\infty}\modd{n}_\rs(t)=\infty) > 0,
    \end{displaymath}
  and hence, from~\eqref{eq:56}, that also 
    \begin{displaymath}
      \Pr(n_\rs(t)\ge\bar{n}_\rs \text{ for all $t\ge0$,}\ %
      \lim_{t\to\infty}n_\rs(t)=\infty) > 0.
    \end{displaymath}
    Hence the process~$\vec{n}(\cdot)$ is transient as required.
\end{proof}

We illustrate the use of the above result with a simple example.

\begin{example}
  \label{ex:6}
  Consider again the network of Example~\ref{ex:2}, in which $R=J=3$.
  As previously observed a necessary and sufficient condition for the
  existence of \emph{some} stable control is given by $2\nu<c$.
  Further, if $3\nu<c$, then Proposition~\ref{prop:1} with the
  Lyapunov function $f$ given by $f(\vec{n})=\sum_{r=1}^3n_r$ shows
  that \emph{any} Pareto efficient control is stable.  Suppose now
  that $2\nu<c$ and that the Pareto efficient control~$\vec{b}$ is
  such that, for $r=1,2$, $b_r(\vec{n})$ is independent of $n_3$ and
  \begin{equation}\label{eq:51}
    b_1(\vec{n}) + b_2(\vec{n}) = c
    \qquad\text{for all $\vec{n}$ such that $\max(n_1,n_2)>0$.}
  \end{equation}
  Thus in particular calls of types~$1$ and $2$ collectively have
  complete priority over calls of type~$3$.  Although we do not in
  this example require any further monotonicity conditions on
  $\vec{b}$, it follows from the requirement of Pareto efficiency that
  the control~$\vec{b}^{\{1,2\}}\colon\Zp^2\to\Rp^3$ is well-defined
  as before, being obtained from $\vec{b}$ by letting $n_3\to\infty$.
  It follows from \eqref{eq:51} that the condition $2\nu<c$
  is necessary and sufficient to ensure that $\vec{b}^{\{1,2\}}$ is
  stable.  We use (a slight modification of)
  Theorem~\ref{thm:priority} to investigate the stability of
  $\vec{b}$.  The stationary distribution $\pi^{\{1,2\}}$ on $\Zp^2$
  induced by $\vec{b}^{\{1,2\}}$ is here just that of the
  process~$(n_1(\cdot),n_2(\cdot))$.  Further, since, from
  \eqref{eq:51}, $n_1(\cdot)+n_2(\cdot)$ is Markov, with a stationary
  distribution which is geometric and independent of any more detailed
  specification of $\vec{b}^{\{1,2\}}$, it follows that
  \begin{equation}
    \label{eq:52}
    \pi^{\{1,2\}}(0,0) = 1 - \frac{2\nu}{c}.
  \end{equation}
  It follows from \eqref{eq:2}, \eqref{eq:51} and the Pareto
  efficiency of $\vec{b}$ that
  \begin{equation}
    \label{eq:53}
      b^{\{1,2\}}_3(0,0)=c,
      \qquad
      b^{\{1,2\}}_3(n_1,0)=0
      \quad\text{for all $n_1\ge1$},
      \qquad
      b^{\{1,2\}}_3(0,n_2)=0
      \quad\text{for all $n_2\ge1$}
  \end{equation}
  We thus have that
  \begin{align}
    \E_{\pi^{\{1,2\}}}b^{\{1,2\}}_3
    & = \sum_{(n_1,n_2)\in\Zp^2}\pi^{\{1,2\}}(n_1,n_2)b^{\{1,2\}}_3(n_1,n_2)
    \nonumber\\
    & \ge c - 2\nu, \label{eq:54}
  \end{align}
  with equality if and only if $b^{\{1,2\}}_3(n_1,n_2)=0$ for all
  $(n_1,n_2)$ such that $\min(n_1,n_2)\ge1$.  But this latter
  condition holds if and only if, for all $(n_1,n_2)$ such that
  $\min(n_1,n_2)\ge1$, we have
  $\min\bigl(b^{\{1,2\}}_1(n_1,n_2),b^{\{1,2\}}_2(n_1,n_2)\bigr)=0$,
  i.e.\ in the case of the control considered in Example~\ref{ex:2} in
  which maximum resource is always allocated to calls of one type,
  and in which we have already observed that we have stability if and
  only if $3\nu<c$.  Otherwise we have strict inequality in
  \eqref{eq:54}.  

  Now note that, although $\vec{b}$ does not here satisfy all the
  conditions for monotonicity given earlier, the assumption that
  $b_1(\vec{n})$ and $b_2(\vec{n})$ are independent of $n_3$ ensures
  that Theorem~\ref{thm:priority} continues to apply, indeed in a
  slightly improved form, to show that the
  condition~$\E_{\pi^{\{1,2\}}}b^{\{1,2\}}_3>\nu$ is necessary and
  sufficient for the stability of $\vec{b}$.  (For the sufficiency,
  note that the proof of part~(i) of the theorem, with $\SS=\{1,2\}$
  and $r'=3$, goes through as before, except that the coupling between
  $\vec{n}(\cdot)$ and process~$\vec{n}^\SS(\cdot)$ is now obtained
  with equality, and so we no longer require the
  condition~\eqref{eq:36} in order to obtain~\eqref{eq:40}.  Similar
  obvious simplifications apply to the proof of part~(ii), which here
  becomes a fairly standard argument and in particular delivers null
  recurrence---and hence instability---in the case
  $\E_{\pi^{\{1,2\}}}b^{\{1,2\}}_3=\nu$.)
  
  Suppose now that $c$ and $\vec{b}$ are held fixed and that $\nu$ is
  allowed to vary.  The obvious coupling argument shows that if
  $\vec{b}$ is stable for any $\nu$ then it is also stable for any
  $\nu'<\nu$.  The above adaptation of Theorem~\ref{thm:priority},
  together with~\eqref{eq:54}, shows that there is some critical
  parameter~$\lambda$ (depending on the detailed specification of
  $\vec{b}^{\{1,2\}}$ and hence $\pi^{\{1,2\}}$) such that
  $1/3\le\lambda\le1/2$ and that $\vec{b}$ is stable if
  $\nu<\lambda{}c$ and unstable if $\nu>\lambda{}c$.  
  For the control of Example~\ref{ex:2} we already know that
  $\lambda=1/3$; otherwise for the case $\nu=c/3$  we have strict
  inequality in \eqref{eq:54} and hence stability; simple continuity
  arguments now give $\lambda>1/3$ in this case.
\end{example}

\section*{Acknowledgement}

The authors are most grateful to Serguei Foss and to Takis
Konstantopoulos for some helpful discussions, and also to the
referee for helpful comments and corrections.


\newpage

\begin{thebibliography}{99}
\bibitem{Ass}
  Asmussen, S. (2003).
  \textit{Applied Probability and Queues.}
  Springer, New York.

\bibitem{BM}
  Bonald, T.\ and Massouli\'e, L. (2001).
  Impact of fairness on Internet performance. In 
  \textit{Proceedings of ACM SIGMETRICS 2001}.
  Cambridge, Massachusetts.

\bibitem{BP}
  Bonald, T.\ and Prouti\`{e}re, A. (2004).
  On Stochastic Bounds for Monotonic Processor Sharing Networks
  \textit{Queueing Systems}, \textbf{47}, 81--106.

\bibitem{BOR}
  Borovkov, A.\ A. (1998).
  \textit{Ergodicity and Stability of Stochastic Processes.}
  Wiley, Chichester.

\bibitem{Dai}
  Dai, J. (1995).
  On positive Harris recurrence of multiclass queueing networks: a
  unified approach via fluid limit models.
  \textit{Annals of Applied Probability}, \textbf{5}, 49--77.

\bibitem{dVKL}
  de Veciana, G., Lee, T.-J.\ and Konstantopoulos, T. (2001).
  Stability and performance analysis of networks supporting elastic
  services.
  \textit{IEEE/ACM Trans.\ on Networking}, \textbf{9}, 2--14.

\bibitem{KMT}
  Kelly, F.\ P., Maulloo, A.\ K.\ and Tan, D.\ K.\ H. (1998).
  Rate control in communication networks: shadow prices, proportional
  fairness and stability.
  \textit{J. Operational Research Society}, \textbf{49}, 237--252.

\bibitem{KW}
  Kelly, F.\ P.\ and Williams, R.\ J. (2004).
  Fluid model for bandwidth sharing.
  \textit{Annals of Applied Probability}, \textbf{14}, 1055--1083.

\bibitem{MW}
  Mo, J.\ and Walrand, J. (2000).
  Fair end-to-end window-based congestion control.
  \textit{IEEE/ACM Trans.\ on Networking}, \textbf{8}, 556--567.

\bibitem{RM} Roberts, J.\ W. \ and Massouli\'e, L (2000).
  Bandwidth sharing and admission control for elastic traffic.
  \textit{Telecommunications Systems}, \textbf{15}, 185--201.

\end{thebibliography}
\end{document}